\documentclass[a4paper,12pt]{article}
\usepackage{amsfonts,amssymb, pdfpages, xcolor}

\title{Asymptotic  study of the initial value problem to a standard one pressure model of  multifluid flows in nondivergence form}
\author{M. Colombeau,\\ \textit{mcolombeau@ime.usp.br}\\
 Instituto de Matem\'atica e Estatistica,\\Universidade  de S\~ao Paulo, Brazil.}

\begin{document}
\maketitle

\begin {abstract}  We construct families of approximate solutions to the initial value problem and provide complete mathematical proofs that they tend to satisfy the standard system of isothermal one pressure two-fluid flows in 1-D when the data are $L^1$ in densities and $L^\infty$ in velocities. To this end, we use a method that reduces this system of PDEs to a family of systems of four ODEs in Banach spaces whose smooth solutions are these approximate solutions. This method is constructive: using standard numerical methods for ODEs one can easily and accurately compute these approximate solutions which, therefore, from the mathematical proof, can serve for comparison with numerical schemes. One observes agreement with previously known solutions from scientific computing [S. Evje, T. Flatten. Hybrid Flux-splitting Schemes for a common two fluid model. J. Comput. Physics 192, 2003, p. 175-210]. We show that one recovers the solutions of these authors (exactly in one case, with a slight difference in another case). Then we propose an efficient numerical scheme for the original system of two-fluid flows and show it gives back exactly the same results as the theoretical solutions obtained above. \\

\end{abstract}
\\
\textit{*this research has been  done thanks to  financial support of FAPESP, processo 2012/15780-9.}\\

\textbf{1.  Introduction}.\\ 
We study a basic model  used to describe mathematically a mixture of two immiscible fluids in the isothermal case and without transfer of momentum between the two fluids,   \cite{EvjeFlatten}  p. 179, \cite{Cortes} p. 465,

\begin{equation}\frac{\partial}{\partial t}(\rho_1 \alpha_1)+\frac{\partial}{\partial x}(\rho_1\alpha_1 u_1)=0,\end{equation}
\begin{equation}\frac{\partial}{\partial t}(\rho_2 \alpha_2)+\frac{\partial}{\partial x}(\rho_2\alpha_2 u_2)=0,\end{equation}
\begin{equation}\frac{\partial}{\partial t}(\rho_1 \alpha_1u_1)+\frac{\partial}{\partial x}(\rho_1\alpha_1( u_1)^2)+\frac{\partial}{\partial x}((p_1-p_1^{int})\alpha_1)+\alpha_1\frac{\partial}{\partial x}(p_1^{int})=g\alpha_1\rho_1,\end{equation}
\begin{equation}\frac{\partial}{\partial t}(\rho_2 \alpha_2u_2)+\frac{\partial}{\partial x}(\rho_2\alpha_2( u_2)^2)+\frac{\partial}{\partial x}((p_2-p_2^{int})\alpha_2)+\alpha_2\frac{\partial}{\partial x}(p_2^{int})=g\alpha_2\rho_2,\end{equation}
\begin{equation}\alpha_1+\alpha_2=1,\end{equation}
\begin{equation}p_1=K_1\rho_1-b_1, \ \ p_2=K_2\rho_2-b_2,\end{equation}
\\
where the two fluids are denoted by the indices $1$  and $2$, for instance mixture of oil and natural gas in extraction tubes of oil exploitation \cite{Avelar}. The physical variables are 
the densities $\rho_i(x,t)$, the velocities $u_i(x,t)$, the volumic proportions
$\alpha_i(x,t)$, the pressures
$p_i(x,t)$, the phasic pressures $p_i^{int}(x,t)$ at the interface, $i=1,2$, and
$g$ is the component of the gravitational acceleration in the direction of the tube.
Equations (6) are the state laws stated  in \cite{EvjeFlatten} p. 179; it is assumed $b_1-b_2>0, K_1>0$ and $ K_2>0$.
Equations (1) and (2) are the continuity equations for each fluid: they express mass conservation. Equations (3) and (4) are the Euler equations for each fluid: they express momentum conservation. A natural assumption is to state the equality of the four pressures $p_i$ and $  p_i^{int}, i=1,2$. This simplest assumption of equal pressure  leads to a nonhyperbolic model, called the equal pressure model  \cite{EvjeFjelde} p.677,  \cite{Munk} p. 2589, \cite{Toumi} p. 287, \cite{Wendroff} p. 372-373 that we study in this paper. \\

We construct families of differentiable functions $S(x,t,\epsilon)$ that, when plugged into the equal pressure model, tend asymptotically to satisfy it when $\epsilon\rightarrow 0$. We prove that these families of functions are weak asymptotic methods. The concept of weak asymptotic method and its relevance has been put in evidence by many authors \cite{Albeverio, Danilov1, Mitrovic, Panov, Shelkovichmat} by explicit calculations and by reduction of the problem of description of nonlinear waves interaction to the resolution of systems of  ordinary differential equations, as a continuation of Maslov' s theory. In other words our families of functions tend to satisfy the system modulo a remainder that tends to 0 when  $\epsilon\rightarrow 0$. To construct these families we use a method which consists in solving a system of four ordinary differential equations in a Banach space whose solutions are the approximate solutions of the one pressure model. This method allows us to compute the solutions with standard convergent numerical schemes for ODEs, thus permitting comparaison with existing numerical solutions of the equal pressure model obtained in scientific computing.  We observe the approximate solutions we obtain agree with the results  presented in \cite{EvjeFlatten}, with a small difference in one case which diminishes in presence of the pressure correction, which can be considered as a mathematical justification of these numerical results. The system (1-6) is in nondivergence form, i.e. the derivatives cannot be transfered to test functions because of the terms $\alpha_i\frac{\partial p_i^{int}}{\partial x} $ in (3, 4). Therefore the study of the solutions of this system in presence of shock waves is  problematic and we use a family of approximate solutions that are classical differentiable functions which permits at the limit to obtain "exact solutions" that are irregular functions such as discontinuous functions. In this way the weak asymptotic methods presented here can be a tool for mathematical and numerical investigations of discontinuous solutions despite this system is in nondivergence form. Various systems in divergence form have been obtained by replacing (10,11) below by their sum, and then by introducing a new equation \cite{Evje1, Evje2, Evje3, EvjeFlattenFriis, EvjeFriis, EvjeKarlsen1, EvjeKarlsen2}.\\

From a physical viewpoint the equations of fluid dynamics are mared with some imprecision since they do not take into account some minor effects and the molecular structure of matter. It is natural to expect these equations and their imprecision should be stated in the sense of distributions in the space variables. Weak asymptotic methods provide approximate solutions that enter into this imprecision for $\epsilon>0$ small enough. Therefore they could be considered as some convenient way to approximate possible solutions to the equations of physics. In the absence of a uniqueness result of a privileged family of weak asymptotic methods (all giving same results) that should represent physics in a given physical situation, we have to content to check numerically that the weak asymptotic methods we present give the known solutions at the limit $\epsilon\rightarrow 0$.\\

\textbf{2. Simplified statement of the system.}\\ In order to simplify the study of the system (1-6) with the equal pressure assumption we transform it into a system of four equations with four unknown functions by changes of unknown and  algebraic calculations. We set
\begin{equation} r_1=\rho_1 \alpha_1, r_2=\rho_2 \alpha_2, \ \alpha=\alpha_1. \end{equation}
Then (1-5) with equal pressures  has the form  

\begin{equation}\frac{\partial}{\partial t}(r_1)+\frac{\partial}{\partial x}(r_1u_1)=0,\end{equation}
\begin{equation}\frac{\partial}{\partial t}(r_2)+\frac{\partial}{\partial x}(r_2u_2)=0\end{equation}
\begin{equation}\frac{\partial}{\partial t}(r_1u_1)+\frac{\partial}{\partial x}(r_1(u_1)^2)+\alpha \frac{\partial}{\partial x}p=gr_1,\end{equation}
\begin{equation}\frac{\partial}{\partial t}(r_2u_2)+\frac{\partial}{\partial x}(r_2(u_2)^2)+(1-\alpha)\frac{\partial}{\partial x}p=gr_2,\end{equation}
and  the two state laws (6) are left unchanged. The  6 unknowns are now $r_1, r_2, u_1, u_2, \alpha$ and $ p$. Then we transform the equations in a way which will be more convenient to construct the weak asymptotic  method since we will have only the four unknown functions $r_1, r_2, u_1$ and $u_2$.\\

$\bullet$ From (6), $p=K_1\rho_1-b_1=K_2\rho_2-b_2 $ implies 
\begin{equation}\rho_2=\frac{-b_1+b_2+K_1\rho_1}{K_2}.\end{equation}
Note that this calculation is linear so it can be done rigorously even in presence of shock waves.\\

$\bullet$\textit{  Calculation of $\alpha$ in function of $r_1$ and $ r_2$}. We multiply the equality $p=K_1\frac{r_1}{\alpha}-b_1=K_2\frac{r_2}{1-\alpha}-b_2$ by $\alpha(1-\alpha)$ to obtain (13) below:
 this is a nonlinear calculation.  In the case of discontinuous solutions it is well known such nonlinear calculations usually change the solutions. This  formal calculation  is usual for this system and the observation of the numerical results in section 6, observation 3, shows a posteriori that this nonlinear calculation giving the formula (13) is justified.
No unjustified nonlinear calculations are done after (13). This calculation gives
\begin{equation}(1-\alpha)K_1r_1-b_1\alpha(1-\alpha)=\alpha K_2r_2-b_2\alpha(1-\alpha),\end{equation}
i.e. $F(\alpha)=0$,
setting 
\begin{equation}F(X)=X^2(b_1-b_2)+X(-K_1r_1-b_1-K_2r_2+b_2)+K_1r_1.\end{equation}
One has $F(0)=K_1r_1> 0$ and $ F(1)=-K_2r_2<0$ which implies that $F$ has one and only one root  $\alpha\in ]0,1[$ in the case $r_1>0$ and $r_2>0$ i.e. in absence of void regions in each fluid. In this case, since $F(1)<0$ and since it is assumed  $b_1-b_2>0$, the second root is $>1$.  Therefore the discriminant $\Delta=(K_1r_1+K_2r_2+b_1-b_2)^2-4(b_1-b_2)K_1r_1$ is $>0$ and the solution $\alpha\in ]0,1[$ is given by 
\begin{equation} \alpha=\frac{K_1r_1+K_2r_2+b_1-b_2-(\Delta)^{\frac{1}{2}}}{2(b_1-b_2)}.\end{equation}

$\bullet$ The following result will be used below
  \begin{equation}0\leq\frac{K_1r_1}{b_1-b_2+K_1r_1+K_2r_2}\leq \alpha\leq 1.\end{equation}
Proof. From (14), $F(X)\geq -X(b_1-b_2+K_1r_1+K_2r_2)+K_1r_1$ since $b_1-b_2>0$; therefore $F(\frac{K_1r_1}{b_1-b_2+K_1r_1+K_2r_2})\geq 0$, hence the result since $F(\alpha)=0$ and $F(1)\leq 0$.$\Box$\\

Finally we can eliminate $\alpha$ from (10, 11) and we obtain the following statement of the system: first the continuity equations
\begin{equation}\frac{\partial}{\partial t}(r_1)+\frac{\partial}{\partial x}(r_1u_1)=0,\end{equation}
\begin{equation}\frac{\partial}{\partial t}(r_2)+\frac{\partial}{\partial x}(r_2u_2)=0,\end{equation}
 then the Euler equations in the form
\begin{equation}\frac{\partial}{\partial t}(r_1u_1)+\frac{\partial}{\partial x}(r_1(u_1)^2)+r_1\frac{\partial}{\partial x}\Phi_1=gr_1,\end{equation}
\begin{equation}\frac{\partial}{\partial t}(r_2u_2)+\frac{\partial}{\partial x}(r_2(u_2)^2)+r_2\frac{\partial}{\partial x}\Phi_2=gr_2,\end{equation}
where 
\begin{equation} \Phi_1=K_1 log \rho_1, \  \rho_1=\frac{r_1}{\alpha}, \  \  \Phi_2=K_2 log \rho_2, \
 \rho_2=\frac{r_2}{1-\alpha}=\frac{-b_1+b_2+K_1\rho_1}{K_2},\end{equation}
with $\alpha$ given by (15).
The system is now a system of four scalar PDEs with the four unknowns $r_1, r_2, u_1$ and $ u_2$.\\

\textbf{3. Statement of the weak asymptotic method.}\\
Setting 
 \begin{equation}u_i^+=\frac{|u_i|+u_i}{2}, \ \ u_i^-=\frac{|u_i|-u_i}{2},\end{equation} 
one has 
\begin{equation} u_i^+-u_i^-=u_i, \ \ u_i^++u_i^-=|u_i|. \end{equation}
The two continuity equations and the two Euler equations are replaced by the following ODEs, $i=1,2$
\begin{equation} \frac{d}{dt}r_i(x,t,\epsilon)=\frac{1}{\epsilon}[(r_i u_i^+)(x-\epsilon,t,\epsilon)-(r_i |u_i|)(x,t,\epsilon)+(r_i u_i^-)(x+\epsilon,t,\epsilon)]+\epsilon^\beta,\end{equation} with $\beta>0$ to be defined later,\\

$ \frac{d}{dt}(r_iu_i)(x,t,\epsilon)=\frac{1}{\epsilon}[(r_i u_iu_i^+)(x-\epsilon,t,\epsilon)-$\begin{equation}(r_iu_i |u_i|)(x,t,\epsilon)+(r_i u_iu_i^-)(x+\epsilon,t,\epsilon)]-r_i(x,t,\epsilon)\frac{\partial}{\partial x}\Phi_i
(x,t,\epsilon)+gr_i(x,t,\epsilon).\end{equation}
The potentials $\Phi_i$,  $i=1,2$, are defined by 

\begin{equation}\Phi_i(x,t,\epsilon)=K_i[log(\rho_i(.,t,\epsilon)+\epsilon^N)*\phi_{\epsilon^\gamma}](x),
\end{equation} with $N$ and $\gamma$ to be defined later, $\phi\in \mathcal{C}_c^\infty(\mathbb{R}), \ \phi\geq 0$ and  $ \int\phi(\mu)d\mu=1$. The convolution in (26) permits that the derivative $\frac{\partial}{\partial x}\Phi_i$ in (25) makes sense: thus the fact the equations (19, 20) are not in divergence form does not cause any trouble for the approximating sequences.
We recall $\alpha$ is defined in (15), \ $\rho_1=\frac{r_1}{\alpha},$ and one will prove $\alpha_i(x,t,\epsilon)>0 \ \forall \epsilon>0$; 
$\rho_2$ is given in (12, 21), $u_i=\frac{r_iu_i}{r_i}$ and one will prove $
r_i(x,t,\epsilon)>0 \ \forall \epsilon>0$. This will follow from (33) below, which, from (13),  implies $\alpha\not=0$ and $\alpha\not=1$.\\

  We assume $r_{i,0}$ and $u_{i,0}, i=1,2$ are given initial conditions on the 1-D torus $\mathbb{T}=\mathbb{R}/(2\pi \mathbb{Z})$  with the properties $r_{i,0} \in L^1(\mathbb{T})$ and  $u_{i,0}\in   L^\infty(\mathbb{T})$ and that $r_{i,0}^\epsilon$ and $ u_{i,0}^\epsilon$ are regularizations of $r_{i,0}$ and $u_{i,0}$ respectively, with  uniform $L^1$ and $L^\infty$ bounds respectively (independent on $\epsilon$), and $r_{i,0}^\epsilon (x)>0 \  \  \forall x$.\\  
  
\textbf{Theorem.}  \textit{If $0<\gamma<\frac{1}{6}$ and $N-1-\beta-3\gamma>0$ the system of four ODEs (24, 25) complemented by the relations (6, 7,  26) provides a weak asymptotic method for the system (17, 18, 19, 20, 21).}\\

 Sections 4 and 5 are devoted to the proof of the Theorem.\\

\textbf{4. A priori inequalities for fixed $\epsilon$.} \\
We  seek solutions  on the 1-D torus $\mathbb{T}=\mathbb{R}/(2\pi \mathbb{Z})$. Families  $(r_{i,0}^\epsilon)_\epsilon$ and  $ (r_{i,0}^\epsilon u_{i,0}^\epsilon)_\epsilon$ of approximations of  initial conditions are given on $\mathbb{T}$. For fixed $\epsilon>0$ we assume existence and uniqueness of a solution of (24, 25) of class $\mathcal{C}^1$  $$[0,\delta(\epsilon)[\longmapsto \mathcal{C}_b(\mathbb{R})^4,$$
$$ \ \ \ \ \ \ \ \   \ \ \ \ \ \ \ \  \ \ \ \ \ \ \ \ t\longmapsto (r_i(.,t,\epsilon),  r_iu_i(.,t,\epsilon))$$
such that
\begin{equation}\exists m>0 \ / \ r_i(x,t,\epsilon)\geq m \ \forall x\in \mathbb{R} \ \forall t\in [0,\delta(\epsilon)[,\end{equation}
\begin{equation}\exists M>0 \ / \| u_i(.,t,\epsilon)\|_\infty\leq M, \| r_i(.,t,\epsilon)\|_\infty\leq M  \ \forall t\in [0,\delta(\epsilon)[.\end{equation}
\\

 \textbf {Proposition 1 (a priori inequalities).}\\
    $ \bullet \  \forall \epsilon>0, \ \forall   t\in [0,\delta(\epsilon)[ \  \ r_i(.,t,\epsilon)\in L^1(\mathbb{T}),$  and   \begin{equation}\int_{-\pi}^\pi  r_i(x,t,\epsilon)dx = \int_{-\pi}^\pi r_i(x,0,\epsilon)dx+2\pi\epsilon^\beta t,\end{equation} 
     \begin{equation}  \bullet  \ \exists C>0 \ / \ \ \|\frac{\partial}{\partial x}\Phi_i(.,t,\epsilon)\|_\infty \leq \frac{C}{\epsilon^{3\gamma}} \ \forall       t\in [0,\delta(\epsilon)[,\forall \epsilon>0, \ \ \ \ \ \ \ \ \ \ \ \ \ \ \ \ \ \ \ \ \ \  \ \  \end{equation}
      \begin{equation} \bullet \|u_i(.,t,\epsilon)\|_\infty \leq  \|u_i(.,0,\epsilon)\|_\infty +\frac{2(C+g)} {\epsilon^{3\gamma}}\delta(\epsilon) \  \forall       t\in [0,\delta(\epsilon)[,\forall \epsilon>0. \ \ \  \end{equation}
      Setting 
  \begin{equation} k(\epsilon)=max_{i=1,2}\|u_i(.,0,\epsilon)\|_\infty+  \frac{2(C+g) \delta(\epsilon)}{\epsilon^{3\gamma}},\end{equation}
   then  $\forall       t\in [0,\delta(\epsilon)[,\forall \epsilon>0, $ 
      \begin{equation} \bullet       \ r_i(x,0,\epsilon)exp(\frac{-k(\epsilon)t}{\epsilon}) \leq r_i(x,t,\epsilon) \leq 2\|r_i(.,0,\epsilon)\|_\infty exp(\frac{2k(\epsilon)t}{\epsilon})  \ \forall x\in \mathbb{R}.\end{equation} 

Proof of  Proposition 1.\\

 $\bullet$ From (23, 24),\\
\\
      $\frac{d}{dt}\int_{-\pi}^{+\pi}r_i(x,t,\epsilon)dx=
      \frac{1}{\epsilon}[\int_{-\pi}^{+\pi}(r_i u_i^+)(x-\epsilon,t,\epsilon)dx-\int_{-\pi}^{+\pi}(r_i u_i^+)(x,t,\epsilon)dx-\\
\\
\int_{-\pi}^{+\pi}(r_i u_i^-)(x,t,\epsilon)dx+\int_{-\pi}^{+\pi}(r_i u_i^-)(x+\epsilon,t,\epsilon)dx]+2\pi\epsilon^\beta 
   =0+2\pi\epsilon^\beta$\\ 
\\
     by periodicity of $r_i$ and $u_i$.  \\

   $\bullet$ From (26),
  $$\frac{\partial}{\partial x}(\Phi_i)(x,t,\epsilon)=K_i\int log[\rho_i(x-y,t,\epsilon)+\epsilon^N]\frac{1}{\epsilon^{2\gamma}}  \phi'(\frac{y}{\epsilon^\gamma})dy.$$
If $\rho_i(x-y,t,\epsilon)\leq 1$, one uses the bound $|log(\epsilon^N)|\leq \frac{const}{\epsilon^\gamma}.$
If $\rho_i(x-y,t,\epsilon)> 1$, one uses the fact that  $\rho_i(.,t,\epsilon)\in L^1(\mathbb{T})$ with $L^1$ norm independent on $\epsilon$ and $t\in[0,\delta(\epsilon)[$. The result that $\rho_i(.,t,\epsilon)\in  L^1(\mathbb{T})$ with $L^1$ norm independent on $\epsilon$ and $t$ follows from formula (16) that implies $\rho_1=\frac{r_1}{\alpha}\leq r_1\frac{b_1-b_2+K_1r_1+K_2r_2}{K_1r_1}$. Then  one notices that  $b_1-b_2>0, K_i>0,r_i>0$ and the result follows from (29). For $\rho_2$ one uses (12). \\ 

 $\bullet$Now we proceed to the proof of (31). 
 From (24) and the assumption that the solution of the ODEs  is of class $\mathcal{C}^1$ on $[0,\delta(\epsilon)[$ valued in the Banach space $\mathcal{C}(\mathbb{T})$, one obtains, for fixed $\epsilon>0$ and for $dt>0$ small enough with $t+dt<\delta(\epsilon)$, that\\

     $r_i(x,t+dt,\epsilon)=r_i(x,t,\epsilon)+$\\
     $$\frac{dt}{\epsilon}[(r_i u_i^+)(x-\epsilon,t,\epsilon)-(r_i |u_i|)(x,t,\epsilon)+(r_i u_i^-)(x+\epsilon,t,\epsilon)]+dt.o(x,t,\epsilon)(dt)+\epsilon^\beta dt=$$ \begin{equation}\frac{dt}{\epsilon}(r_i u_i^+)(x-\epsilon,t,\epsilon)+(1-\frac{dt}{\epsilon}|u_i|(x,t,\epsilon))r_i(x,t,\epsilon)+\frac{dt}{\epsilon}(r_i u_i^-)(x+\epsilon,t,\epsilon)+dt.o(x,t,\epsilon)(dt)+\epsilon^\beta dt\end{equation}
     \\ 
      where $\|o(.,t,\epsilon)(dt)\|_\infty \rightarrow 0$ when $dt\rightarrow 0$ uniformly for $t$ in a compact set of $[0,\delta(\epsilon)[$, from the mean value theorem in the form $f(t+dt)=f(t)+dt f'(t) +dt.  r(t,dt),$ with $ \|r(t,dt)\|\leq sup_{0<\theta<1}\|f'(t+\theta dt)-f'(t)\|$. Notice that there is no uniformness in $\epsilon$. For $dt>0$ small enough (depending on $\epsilon$)   the single term $  (1-\frac{dt}{\epsilon}|u_i|(x,t,\epsilon))r_i(x,t,\epsilon)$ dominates the term $dt.o(x,t,\epsilon)(dt)$ from (27, 28). Since, further, $r_i u_i^\pm\geq0$, one can invert (34). Dropping  the useless term $\epsilon^\beta dt$ one obtains\\
 
$\frac{1}{r_i(x,t+dt,\epsilon)}\leq$\\
\\
$[\frac{dt}{\epsilon}(r_i u_i^+)(x-\epsilon,t,\epsilon)+[1-\frac{dt}{\epsilon}|u_i|(x,t,\epsilon)]r_i(x,t,\epsilon)+\frac{dt}{\epsilon}(r_i u_i^-)(x+\epsilon,t,\epsilon)]^{-1}+\\
\\
dt.o(x,t,\epsilon)(dt)$\\
\\
where the new $o$ has still the property that $\|o(.,t,\epsilon)(dt)\|_\infty \rightarrow 0$ when $dt\rightarrow 0$ uniformly for $t\in[0,\delta'] $  if $\delta'<\delta(\epsilon)$.\\

      Applying the analog of (34) for $r_i u_i$ in place of $r_i$, with the supplementary terms $r_i\frac{\partial}{\partial x}(\Phi_i)$ and $gr_i$ from (25),    one obtains, using (27, 28)
 $$u_i(x,t+dt,\epsilon)=\frac{(r_i u_i)(x,t+dt,\epsilon)}{r_i(x,t+dt,\epsilon)}\leq$$
$$ \frac
 {\frac{dt}{\epsilon}(r_i u_i u_i^+)(x-\epsilon,t,\epsilon)+[1-\frac{dt}{\epsilon}|u_i|(x,t,\epsilon)](r_i u_i)(x,t,\epsilon)+\frac{dt}{\epsilon}(r_i u_i u_i^-)(x+\epsilon,t,\epsilon)}
 {\frac{dt}{\epsilon}(r_i u_i^+)(x-\epsilon,t,\epsilon)+[1-\frac{dt}{\epsilon}|u_i|(x,t,\epsilon)]r_i(x,t,\epsilon)+\frac{dt}{\epsilon}(r_i u_i^-)(x+\epsilon,t,\epsilon)}$$
 \begin{equation} +dt\frac{r_i(x,t,\epsilon)}{r_i(x,t+dt,\epsilon)}[|\frac{\partial}{\partial x}(\Phi_i)(x,t,\epsilon)|+g]+ dt.o(x,t,\epsilon)(dt)\end{equation}
where the new $o$ has the same property as in (34)  for fixed $\epsilon$.
 For $dt>0$ small enough the first term  in the second member is  a  barycentric combination of $u_i(x-\epsilon,t,\epsilon),  u_i(x,t,\epsilon)$ and $ u_i(x+\epsilon,t,\epsilon)$. The quotient  $\frac{r_i(x,t+dt,\epsilon)}{r_i(x,t,\epsilon)} $ tends to $1$ when $dt\rightarrow 0$ (for fixed $\epsilon$). Finally one obtains, using also (30),  that 
 \begin{equation} \|u_i(.,t+dt,\epsilon)\|_\infty\leq  \|u_i(.,t,\epsilon)\|_\infty+dt\frac{const}{\epsilon^{3\gamma}}+dt.\|o(.,t,\epsilon)(dt)\|_\infty\end{equation}
 with uniform bound of $o$ when $t$ ranges in a compact set in $[0,\delta(\epsilon)[$,  for fixed $\epsilon$. One obtains the bound (31) as in \cite{Colombeaugravitation} by dividing the interval $[0,t]$ into $n$ small intervals $[\frac{it}{n},\frac{(i+1)t}{n}], 0\leq i \leq n-1$, applying (36) in each small interval, which gives
 $$\|u_i(.,(i+1)\frac{t}{n},\epsilon)\|_\infty \leq \|u_i(.,i\frac{t}{n},\epsilon)\|_\infty+\frac{t}{n}\frac{const}{\epsilon^{3\gamma}}+\frac{t}{n}o(\frac{t}{n}),$$
 summing on $i$ and using that $o(\frac{t}{n})\rightarrow 0$ when $n\rightarrow \infty$.\\
 
 $\bullet$The proofs of the two inequalities (33) follows from (24) that gives the inequalities $\frac{d}{dt}r_i(x,t,\epsilon)\geq-\frac{\|u_i\|_\infty}{\epsilon}r_i(x,t,\epsilon)$ and $\frac{d}{dt}r_i(x,t,\epsilon)\leq \frac{2\|u_i\|_\infty\|r_i\|_\infty}{\epsilon}$ using (31) to evaluate $\|u_i\|_\infty$. They are  given in detail in section 2 of \cite{Colombeaugravitation}. \\

 The existence of a unique global solution to (24, 25) for fixed $\epsilon$ is obtained from the a priori estimates (29-33)  from classical ODEs arguments of the theory of ODEs in Banach spaces in the Lipschitz case. Indeed for fixed $\epsilon>0$ , if $0<\lambda<1$ and $\Omega_\lambda:=\{(X_i,Y_i)\in \mathcal{C}(\mathbb{T})^4 /  \ \forall x\in \mathbb{T} \ \lambda<X_i(x)<\frac{1}{\lambda}, |Y_i(x)|<\frac{1}{\lambda}\}$ the four equations (24, 25) with variables $X_i=r_i$ and $ Y_i=r_i u_i$ have the Lipschitz property on $\Omega_\lambda$ with values in $\mathcal{C}(\mathbb{T})^4$, with Lipschitz constants $\leq \frac{1}{\lambda^3}$. We refer to \cite {Colombeaugravitation} section 4 for details.\\

\textbf{5. Proof of the weak asymptotic method.} \\
  It remains to prove that the solution of the system of ODEs (24, 25) and the formula (26)  provide a weak asymptotic method for system (17, 18, 19, 20, 21) when $\epsilon\rightarrow 0$. To this end one has to prove that $\forall \psi\in \mathcal{C}_c^\infty(\mathbb{R})$, (37-39) below hold when $\epsilon\rightarrow 0$
 \begin{equation} \int \frac{d}{dt}r_i(x,t,\epsilon)\psi(x)dx-\int (r_i u_i)(x,t,\epsilon)\psi'(x)dx \rightarrow 0,\end{equation}
  $ \int \frac{d}{dt}(r_i u_i)(x,t,\epsilon)\psi(x)dx-$\begin{equation}\int (r_i (u_i)^2)(x,t,\epsilon)\psi'(x)dx +\int r_i(x,t,\epsilon)\frac{\partial}{\partial x}(\Phi_i)(x,t,\epsilon)\psi(x)dx-g\int r_i(x,t,\epsilon)\psi(x)dx\rightarrow 0,\end{equation}
  \begin{equation}  \int\Phi_i(x,t,\epsilon)\psi(x)dx-K_i\int  \log[\rho_i(x,t,\epsilon)]\psi(x)dx\rightarrow 0 \end{equation} 
 where (37) means satisfaction of (17, 18), (38) satisfaction of (19, 20) and (39) satisfaction of the two state laws in (21) in the sense of distributions at the limit $\epsilon\rightarrow 0$.\\ 

  The proof of (37) is as follows: from (23, 24, 29, 30, 31), a change of variable and $\frac{\psi(x+\epsilon)-\psi(x)}{\epsilon}=\psi'(x)+O_x(\epsilon)$,\\
     
  $ \int \frac{d}{dt}r_i(x,t,\epsilon)\psi(x)dx=
  \frac{1}{\epsilon}\int (r_i u_i^+)(x,t,\epsilon)[\psi(x+\epsilon)-\psi(x)]dx -\frac{1}{\epsilon}\int (r_i u_i^-)(x,t,\epsilon)$\\
\\
$[\psi(x)-\psi(x-\epsilon)]dx+  
\int\epsilon^\beta \psi(x)dx =
  \int (r_i u_i)(x,t,\epsilon)\psi'(x)dx+\int_{compact} (r_i u_i^+)(x,t,\epsilon)$ \\
\\
$O_x(\epsilon)dx+ 
\int_{compact} (r_i u_i^-)(x,t,\epsilon) O_x(\epsilon)dx+O(\epsilon^\beta)= \int(r_iu_i)(x,t,\epsilon)\psi'(x) dx+$\\
\\
$(const+\frac{const}{\epsilon^{3\gamma}}) O(\epsilon)+O(\epsilon^\beta)= \int(\rho u)(x,t,\epsilon)\psi'(x) dx+O(\epsilon^{1-3\gamma})+O(\epsilon^\beta)$.\\
  \\
   This gives (37) if $0<\gamma<\frac{1}{3}$. The proof of (38) is similar since the additional terms $\int r_i(x,t,\epsilon)\frac{\partial}{\partial x}\Phi_i(x,t,\epsilon)\psi(x)dx$ and $g\int r_i(x,t,\epsilon)\psi(x) dx$ are the same in (25) and (38): one obtains a remainder $\frac{const}{\epsilon^{6\gamma}}  O(\epsilon)$ because of one more factor $u_i$ and the bound (31). Finally one chooses $0<\gamma<\frac{1}{6}$.\\


 To check (39) one has to prove from (26) that $\forall \psi \in \mathcal{C}_c^\infty(\mathbb{R})$
  \begin{equation}\int\{[(log(\rho_i(.,t,\epsilon)+\epsilon^N)*\phi_{\epsilon^\gamma}](x)-log[\rho_i(x,t,\epsilon)]\}\psi(x)dx \rightarrow 0   \end{equation}
  when $\epsilon\rightarrow 0$. To this end we share the integral (40) into two parts (41, 42) below and we prove that each  tends to 0 when $\epsilon\rightarrow 0$. Let 
   \begin{equation}I=  \int\{[(log(\rho_i(.,t,\epsilon)+\epsilon^N)*\phi_{\epsilon^\gamma}](x)-log[\rho_i(x,t,\epsilon)+\epsilon^N]\}\psi(x)dx \end{equation} 
  and 
  \begin{equation} J=\int\{(log[\rho_i(x,t,\epsilon)+\epsilon^N]-log[\rho_i(x,t,\epsilon)]\}\psi(x)dx.\end{equation}
  Now\\

  $I=\int\{(log[\rho_i(x-\epsilon^\gamma\mu,t,\epsilon)+\epsilon^N]-log[\rho_i(x,t,\epsilon)+\epsilon^N]\}\phi(\mu)\psi(x)d\mu dx=$\\
\\
  $\int log[\rho_i(x,t,\epsilon)+\epsilon^N]\phi(\mu)[\psi(x+\epsilon^\gamma\mu)-\psi(x)]d\mu dx.$\\
  \\
  Since $\rho_i(x,t,\epsilon)\geq 0$ from (16, 21, 33), using its $L^1$ property (29) in the case $\rho_i(x,t,\epsilon)>1$ and using the term $\epsilon^N$ in the case  $\rho_i(x,t,\epsilon)\leq1$, as in the proof of (30),  one has  $|I|\leq const.log(\frac{1}{\epsilon})\epsilon^\gamma$. Therefore $I\rightarrow 0$ when $\epsilon\rightarrow 0$.\\
  
  Now, (42) and the mean value theorem give 
  \begin{equation}|J|\leq \epsilon^N\frac{1}{min (\rho_i)} const\end{equation}
  if $min(\rho_i)$ denotes the inf of $\rho_i(x,t,\epsilon)$ for fixed  $t,\epsilon$ when $x$ ranges in a compact set containing the support of $\psi$.
  The problem is to obtain an inf. bound of $min(\rho_i)$; this is the purpose of the term $\epsilon^\beta$ in (24).  From  (24),
  $\frac{d r_i}{dt}(x,t,\epsilon)\geq-\frac{1}{\epsilon}r_i(x,t,\epsilon)\|u_i(.,t,\epsilon)\|_\infty +\epsilon^\beta \geq - \frac{1}{\epsilon}r_i(x,t,\epsilon)\frac{const}{\epsilon^{3\gamma}}T +\epsilon^\beta$ if $t\in [0,T[$, from (31) applied with $\delta(\epsilon)=T$.\\

  Setting $A:=const \frac{T}{\epsilon^{1+3\gamma}}$ and   $B:=\epsilon^\beta$, one has 
  $\frac{dr_i}{dt}\geq -Ar_i+B. $ Comparing with the exact solution of 
    the ODE $\frac{dX}{dt}(x,t)=-AX(x,t)+B$ with initial condition $X(x,0)=r_{i,0}(x,\epsilon)$,  namely
   $X(x,t)=r_{i,0}(x,\epsilon) e^{-At}+\frac{B}{A}(1-e^{-At})\geq \frac{B}{A}(1-e^{-At})$, we obtain the bound
  
  \begin{equation}r_i(x,t,\epsilon)\geq const(t).\epsilon^{1+\beta+3\gamma}\end{equation}
for $\epsilon>0$ small enough and fixed $t$.
Now using (16) we can obtain a lower bound of $min\rho_i$\\

 \begin{equation}\rho_1(x,t,\epsilon) =\frac{r_1(x,t,\epsilon)}{\alpha(x,t,\epsilon)}\geq r_1(x,t,\epsilon)\geq const. \epsilon^{1+\beta+3\gamma}. \end{equation}
Similarly, from (12)\begin{equation}\rho_2(x,t,\epsilon)  \geq const.\epsilon^{1+\beta+3\gamma}.\end{equation}

  From (43),
  $|J|\leq const(t).\epsilon^{N-1-\beta-3\gamma}$  and  it suffices to choose $N-1-\beta-3\gamma>0$ 
  to obtain that $J\rightarrow 0$ when $\epsilon \rightarrow 0.$ $\Box$\\


\textbf{6. Numerical observations from the weak asymptotic method.} \\
We will present two shock tube problems selected from \cite{EvjeFlatten}. The pressure laws are those in \cite{EvjeFlatten} p. 179-180: $K_1=10^6, K_2=10^5, b_i=K_i\rho_{0,i}-p_{0,i}, \rho_{0,1}=1000, p_{0,1}=10^5, \rho_{0,2}=0$ and  $p_{0,2}=0$. The final time is $T=0.001$, with 1000 cells on [0,1] (or $T=0.01$ on $[0,10]$), therefore $\Delta x=\epsilon=(1000)^{-1} $ and the CFL number is $r=\frac{\Delta t}{\Delta x}= 10^{-6}$. We use the explicit Euler order one method for the ODEs (24, 25). We choose  $\delta=1, \beta=100$ and $ N=100$ ($\beta$ and $N$ do not matter since there is no void region in any fluid). We regularize the initial conditions $\omega_{i,0}=r_{i,0}, r_{i,0}u_{i,0}$ by an averaging

 \begin{equation}\nu \omega_{i,0}(x-\epsilon)+(1-2\nu)\omega_{i,0}(x)+\nu\omega_{i,0}(x+\epsilon), \
 \nu=0.1.\end{equation}  
We represent the convolution   in (26) by a similar averaging of  $\Phi_i$ on 5 cells instead of 3 with coefficient $\nu$=0.15. We use a small averaging as (47) at each step in $r_i$ and $r_i u_i, i=1,2$. Concerning this last averaging one observes that the minimal needed values of $\nu$ tend to 0 when $r\rightarrow 0$: $\nu=10^{-2}, 10^{-3}$ and $ 10^{-4}$ when $r=10^{-4}, 10^{-5}$ and $ 10^{-6}$ respectively. Therefore this regularization can be considered as a numerical artefact absent in the ODE formulation which corresponds to $r=0$. 
The Riemann conditions are $\alpha=0.71,0.7,    p=265000,265000,   u_1=1,1$ and  $u_2=65,50$  for test 1 and $\alpha=0.7,0.1,    p=265000,265000,   u_1=10,15$ and  $ u_2=65,50$  for test 2.\\

\textit{Observation 1}. For shock tube problem 1 (figure 1) one observes  the same results as those depicted in $\cite{EvjeFlatten}$. For shock tube problem 2 (figure 2) one observes a slight difference for the second step  value in the right panels: 2.46 $10^5$ instead of 2.50 $10^5$ (top panel) and 89 instead of 84 (bottom panel). These values do not change with discretizations ranging from 100 to 20000 cells, with different values of $r$ and the other parameters, and are also  exactly those obtained from the direct adaptation of the scheme in  section 7 below. With the pressure correction adopted in \cite{EvjeFlatten} one observes from the scheme in section 7 that this difference tends to disappear, figures 4 and 5, therefore it is presumably a consequence of  the pressure correction adopted in \cite{EvjeFlatten}. Modulo this difference one observes  that the results we obtain without pressure correction agree with the results obtained in \cite{EvjeFlatten} even with pressure correction, which appear therefore as depictions of approximate solutions. \\

\textit{Observation 2}. In both tests it has been observed that the left and right discontinuities satisfy with great precision the 3 standard jump conditions of system (8-11): the two ones from (8, 9) and the third one from the equation obtained by adding (10) and (11). They satisfy also with great precision the two formal jump conditions (62) one can calculate from nonlinear algebraic calculations with the nonconservative equations as done in the appendix.  

The arrays below give on a line the values of the wave velocities computed from the two equations (8, 9) i.e. $c=\frac{[r_1u_1]}{[r_1]}$ and  $c=\frac{[r_2u_2]}{[r_2]},$  the value computed  adding the equations (10, 11) without gravitation, i.e. $c=\frac{[r_1(u_1)^2+r_2(u_2)^2+p]}{[r_1u_1+r_2u_2]}$, and the two formal results (62). They are calculated from the numerical step values in figures 1 and 2. We first give the results for the shock tube problem 1, then for the shock tube problem 2. When the jump conditions are satisfied all values on a line should be equal since they are the value of the velocity of a shock wave obtained from the 5 different formulas.

\vskip13cm

 \begin{figure}[h]
\centering
\includepdf[width=\textwidth]{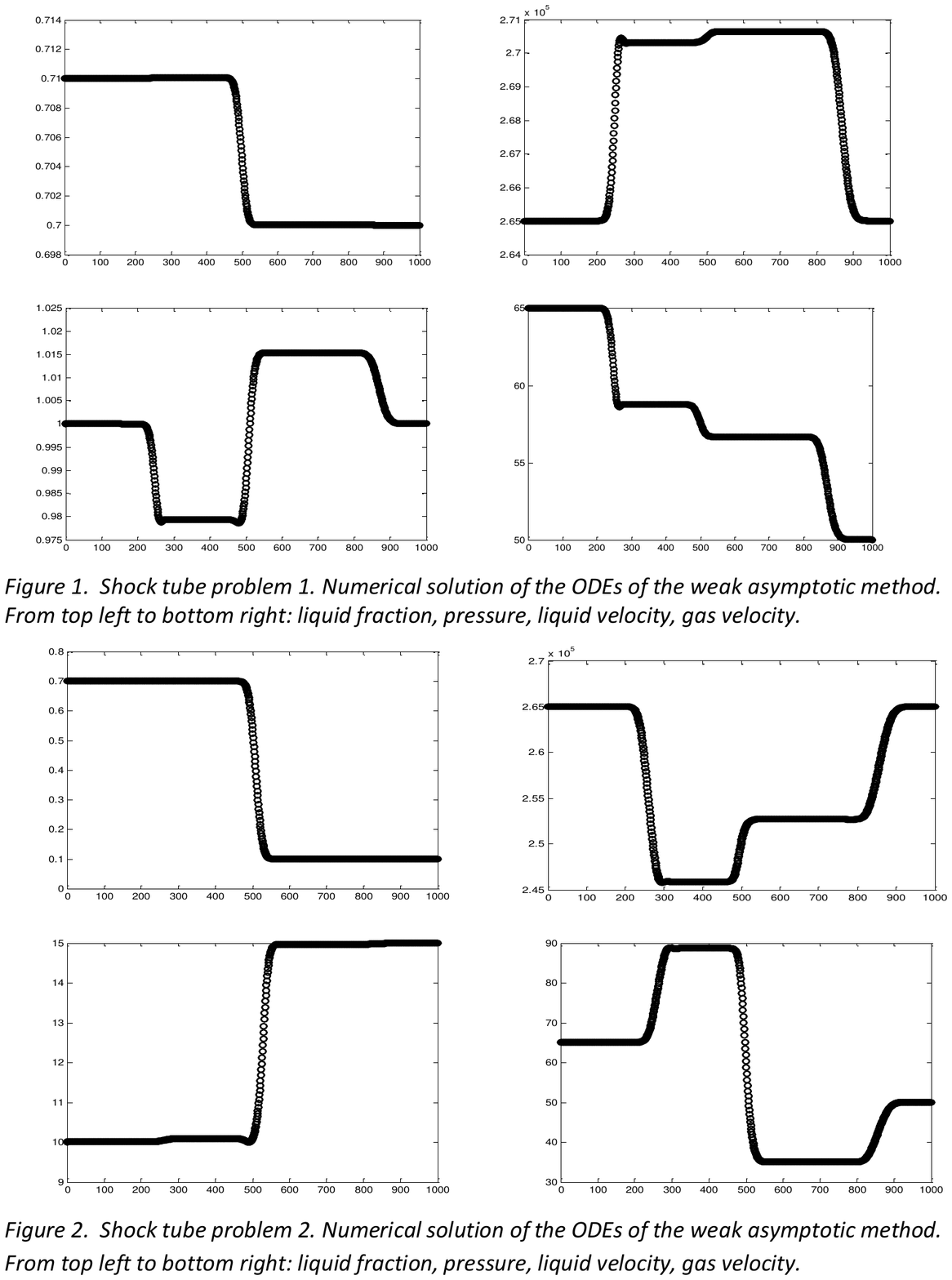}
\end{figure}

\begin{displaymath}
\begin{tabular}{|c||c|c|c|c|c|}
\hline
\multicolumn{6}{|c|}{shock tube problem 1 }\\
\hline 
     &$ c_1$&$c_2$&$c_3$&$c_4$&$c_5$\\ 
 \hline  
left&-255.95 & -255.90 & -255.82 &-255.84  &-255.75\\ 
 \hline
middle& -1.52&-1.49&17.94&10.23&-0.95\\
 \hline
 right&370.35&370.23& 370.24& 370.28&370.15\\
 \hline
 \end{tabular} 
\end{displaymath}

\begin{displaymath}
\begin{tabular}{|c||c|c|c|c|c|}
\hline
\multicolumn{6}{|c|}{shock tube problem 2}\\
\hline 
     &$ c_1$&$c_2$&$c_3$&$c_4$&$c_5$\\   
 \hline
left& -240.23 & -240.85 & -240.81 &-240.65  &-240.91\\
 \hline
middle& 9.34&9.30&8.30&13.78&10.52\\
 \hline
 right&358.32&358.70& 358.84& 358.53&358.77\\
 \hline

 \end{tabular} 
\end{displaymath}

 These jump formulas  are very well satisfied by the left and right discontinuities but are not  satisfied by the middle discontinuity except the two jump conditions from the two continuity equations. Since it is proved the results depict approximate solutions (from the theorem and from a careful numerical solution of the ODEs) an explanation could be that the middle discontinuity is not a classical shock wave as suggested by the   singularities  often observed on top or bottom of this discontinuity,  which could denote it is some kind of  more complicated wave, possibly not a  shock wave. To test this hypothesis we did numerical tests for different  volumic compositions of the fluids. One observes that in the case of equal volume fractions on both sides in the Riemann problem there appears a very neat singularity in the middle discontinuity which is present in the other cases but far clearly visible on the volumic fraction  when both sides of the volumic fraction are equal (figure 3). With the values of pressure and velocities of shock tube problem 1 the observed singularity in volumic fraction is small, while it is quite large with the values of shock wave problem 2, figure 3. This explains why the middle discontinuity does not satisfy well, or does not satisfy at all  in some cases, the expected conservative jump conditions: it is not a classical shock wave i.e. a mere moving discontinuity. It is natural that something else than a mere shock wave occurs: if the Riemann problem were solved by three standard shock waves we would have 8+3=11 unknown values (the 3 velocities and the 8 step values) for 12 equations (the 4 jump conditions at each discontinuity supposing one has solved the ambiguity in the 2 nonconservative equations).  The values of wave velocities in case of figure 3, computed from the 5 algebraic formulas as in the above arrays  are

\begin{displaymath}
\begin{tabular}{|c||c|c|c|c|c|}
\hline
\multicolumn{6}{|c|}{shock tube problem 3}\\
\hline 
     &$ c_1$&$c_2$&$c_3$&$c_4$&$c_5$\\   
 \hline
left& -253.35 & -253.34 & -253.33 &-253.33  &-253.31\\
 \hline
middle& -8473&-1099&25.3&12.8&-16.9\\
 \hline
right& 368.99&368.96& 368.96& 368.97&368.92\\
 \hline

 \end{tabular} \end{displaymath}
 \\     
 Since we have an approximate solution that can be computed with arbitrary precision it is possible to observe this singular part of the solution. Numerical investigation on the "object" that appears in the liquid fraction for $\alpha=0.60$ (top-left panel in figure 3) shows that \\
\\
\\
\\
\\
\\
\\
\\
\\
\\
\\
\\
\\
\\
\\
\\
\\
\\
\\
\\
\\
\\
$\bullet$ The area of the region located below the line $\alpha=0.60$ and above the curve $\alpha$ is  constant (independent on $\epsilon$) for fixed time when $\epsilon$ varies and it is proportional to time even up to very large values of time (tests were done up to 100 times the value of time used in figures 1, 2 and 3 with the scheme in  section 7).

$\bullet$ For rather small values of the time such as those in figures 1, 2 and 3 the object travels with constant speed and  its width on $\alpha=0.60$ tends to 0 when $\epsilon\rightarrow 0$, roughly as  $\sqrt\epsilon$ for fixed $t$ and as $\sqrt t$ for fixed $\epsilon$;  its minimum value diminishes when the time increases. For large values of the time this decrease of the minimum is stopped because one always has $\alpha(x)>0 \  \forall x$ and one observes the width of the object then increases proportionally to $t$ so as to maintain an area proportional to time.

 \begin{figure}[h]\centering\includepdf[width=\textwidth]{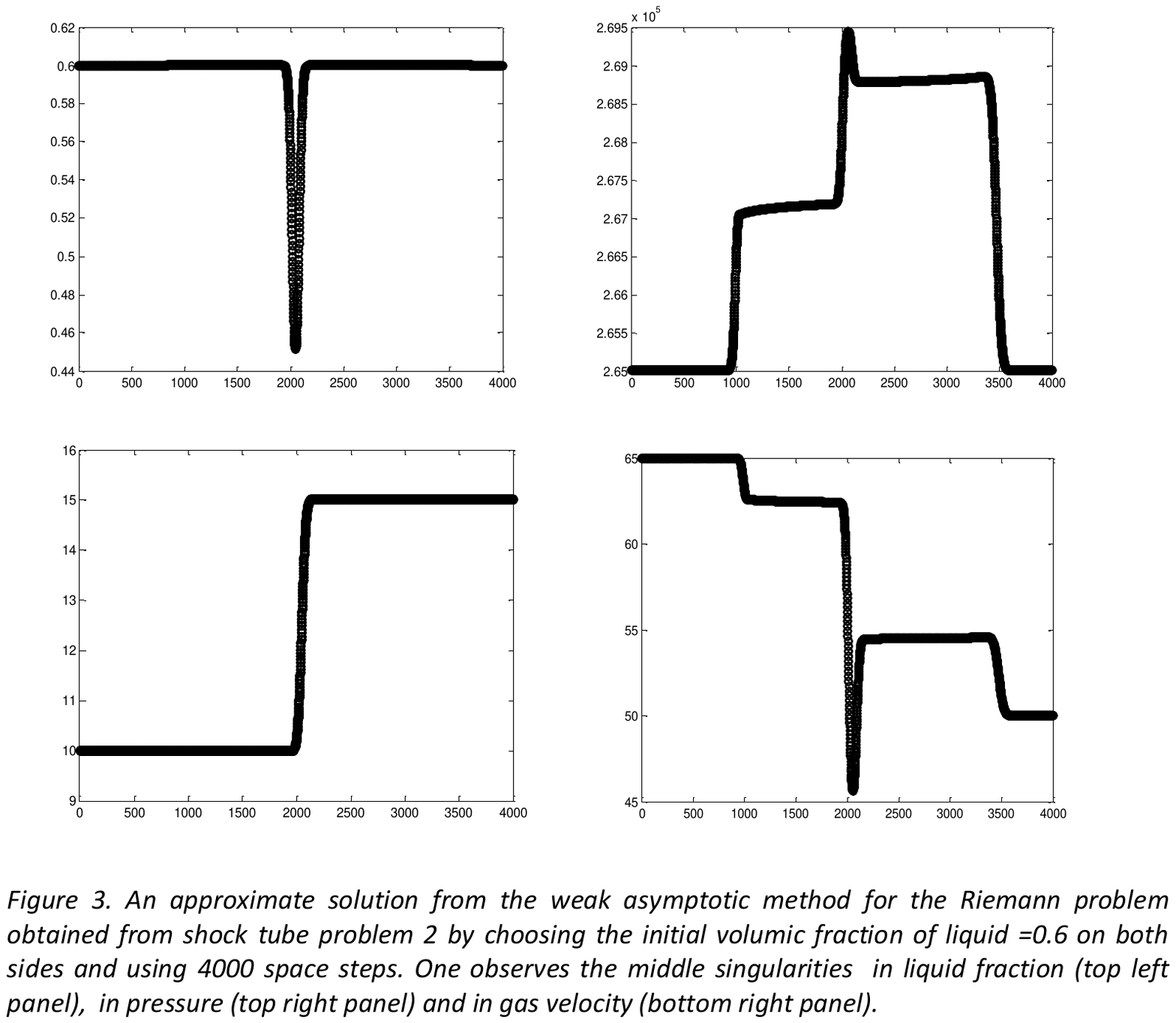}\end{figure}

\textit{Observation 3.} One could state different values of $\gamma$ for the spreading of the state laws of the two fluids when they are very different, for instance a liquid and a gas: we have observed that modifications representing the convolution are unefficient to produce significative differences in the solution because in the present case the discontinuities take place on a large number of cells thus making the results rather unsensitive to modifications that would not be important enough to modify significatively  the aspect of the jumps of $\alpha$ and $p$ (the nonconservative terms in (10, 11)): indeed the great sensibility on the slight modifications of the  schemes for systems in nondivergence form  has been observed  in the case the discontinuities take place on a very small number of cells. The numerical schemes observed  in the case of the multifluid system are robust in the sense that small modifications of the scheme do not affect significantly the result precisely because the discontinuities are spread over a large number of cells. Because of this fact one observes in the three arrays corresponding to figures 1, 2 and 3 that not only the conservative jump conditions (the three  values $c_1,c_2$ and $c_3$) but also the two formal jump conditions (the two values $c_4$ and $c_5$) are satisfied, showing the evidence that, to some extent, one can compute formally on the system, thus allowing the formal nonlinear calculation done to obtain formula (13).\\

\textbf{7. A  transport-correction scheme.}\\
We propose here a natural numerical scheme for the numerical solution of system (17-21).   The scheme is  an adaptation of the  Le Roux et al numerical method of splitting  into transport  and pressure correction  as described in \cite{Baraille},  extending to two fluids the scheme done in \cite{ColombeauNMPDE} for one  fluid. The space $\mathbb{R}\times [0,+\infty[$  is divided into rectangular cells $[ih-\frac{h}{2},ih+\frac{h}{2}]\times [n\Delta t, (n+1)\Delta t[, i\in\mathbb{Z}, n\in\mathbb{N}$.\\

Given the family $\{(r_1)_i^n,(r_2)_i^n,(r_1u_1)_i^n,(r_2u_2)_i^n\}_{i\in\mathbb{Z}}$ of values of these variables on the interval $[ih-\frac{h}{2},ih+\frac{h}{2}]$ at time $n\Delta t$ we seek the family   of values $\{(r_1)_i^{n+1},(r_2)_i^{n+1},(r_1u_1)_i^{n+1},(r_2u_2)_i^{n+1}\}_{i\in\mathbb{Z}}$ at time $(n+1)\Delta t$.\\

$\bullet$ \textit{First step: transport}. For k=1, 2
  
\begin{equation}(u_k)_i^n:=\frac{(r_ku_k)_i^n}{(r_k)_i^n} \end{equation}  if  $ (r_k)_i^n\not=0$, any value if  $(r_k)_i^n=0,$
\begin{equation} (u_k)_i^{n,+}:=\frac{|(u_k)_i^n|+(u_k)_i^n}{2},  (u_k)_i^{n,-}:=\frac{|(u_k)_i^n|-(u_k)_i^n}{2}.\end{equation} 
The CFL condition is $r|(u_k)_i^n|<1 \   \forall k,i,n$.
Then if  $r=\frac{\Delta t}{h}$ we set
\begin{equation} (\overline{r_k})_i:=r (r_k)_{i-1}^n(u_k)_{i-1}^{n,+}+(1 -r |(u_k)_i^n|) (r_k)_i^n  + r (r_k)_{i+1}^n(u_k)_{i+1}^{n,-},\end{equation} 
\begin{equation} (\overline{r_ku_k})_i:=r (r_ku_k)_{i-1}^n(u_k)_{i-1}^{n,+} +(1 -r |(u_k)_i^n|) (r_ku_k)_i^n   + r (r_ku_k)_{i+1}^n(u_k)_{i+1}^{n,-}.\end{equation}\\

 $\bullet$ \textit{Second step: averaging}. We choose a value $\mu, 0<\mu<0.5$,

\begin{equation} (r_k)_i^{n+1}:=\mu(\overline{r_k})_{i-1}+(1-2\mu)(\overline{r_k})_{i}+\mu(\overline{r_k})_{i+1},\end{equation}

\begin{equation}\widetilde{ (r_ku_k)}_i:=\mu(\overline{r_ku_k})_{i-1}+(1-2\mu)(\overline{r_ku_k})_{i}+\mu(\overline{r_ku_k})_{i+1}.\end{equation}\\

$\bullet$ \textit{Third step: pressure correction}. 
\begin{equation} \Delta_i:=(K_1\overline{ (r_1)_i}+K_2\overline{ (r_2)_i}+b_1-b_2)^2-4(b_1-b_2)K_1\overline{ (r_1)_i},\end{equation}
\begin{equation} \alpha_i:=\frac{K_1\overline{ (r_1)_i}+K_2\overline{ (r_2)_i}+b_1-b_2-\sqrt(\Delta)}{2(b_1-b_2)},\end{equation}

\begin{equation} p_i:=K_1\frac{  \overline{(r_1)_i}  } {  \alpha_i  }-b_1 \ \  if \  \alpha_i\not=0,\end{equation}
\\
\begin{equation} (r_1u_1)_i^{n+1}:=\widetilde{ (r_1u_1)}_i-\frac{r}{2} \alpha_i(p_{i+1}-p_{i-1}),\end{equation}
\\
\begin{equation} (r_2u_2)_i^{n+1}:=\widetilde{ (r_2u_2)}_i-\frac{r}{2} (1-\alpha_i)(p_{i+1}-p_{i-1}).\end{equation}
\  \  \\
In (52, 58) we have obtained  the family   $\{(r_1)_i^{n+1},(r_2)_i^{n+1},(r_1u_1)_i^{n+1},(r_2u_2)_i^{n+1}\}_{i\in\mathbb{Z}}$.\\
 \\

Now we justify the choice of an arbitrary value in density when a denominator in (48) is null.\\

\textbf{Proposition.} 
\textit{When $(r_k)_i^{n+1}=0,  k=1$ or $2$,  then  $(r_ku_k)_i^{n+1}=0$. }\\

proof.
Assume $(r_k)_i^{n+1}=0$.   Then from (52), the strict inequality in $\mu$ and the positiveness of $r_k$ imply
 \begin{equation}(\overline{r_k})_{i-1}=0=(\overline{r_k})_{i}=(\overline{r_k})_{i+1}.\end{equation} 
 Now notice that $(\overline{ r_k})_i=0$ implies $(r_k)_i^n=0$ from (50) and the strict inequality in the CFL condition.  Further since $r\not=0$ it also implies from (50) that $ (r_k)_{i-1}^n(u_k)_{i-1}^{n,+}=0$, which implies $ (r_k u_k)_{i-1}^n(u_k)_{i-1}^{n,+}=0$, and  similarly
$ (r_k u_k)_{i+1}^n(u_k)_{i+1}^{n,-}=0$. Therefore $(\overline{r_k})_i=0$ implies  $(\overline{r_ku_k})_i=0$. Therefore from (59) 
\vskip10cm
\textit{Figure 4. The shock tube problem 2 without correction (continuous curve) and with correction (+). One observes a small difference in two step values in the right panels.}
 \begin{equation}(\overline{r_ku_k})_{i-1}=0=(\overline{r_ku_k})_{i}=(\overline{r_ku_k})_{i+1}. \end{equation}  
Therefore  from (53) $ (\widetilde{r_ku_k})_i=0$.  From (56, 59) one has also $p_{i-1}=b_1=p_{i}=p_{i+1}$ if $k=1$. Finally, from (57, 58), we obtain $(r_ku_k)_i^{n+1}=0$. $\Box$\\

Following calculations in \cite{ColombeauSiam, ColombeauNMPDE} one can prove, under assumptions to be checked, such as boundedness of the velocity field when $h\rightarrow 0$, that the scheme tends to satisfy the equations when $h\rightarrow 0$.
 \begin{figure}[h]\centering\includepdf[width=\textwidth]{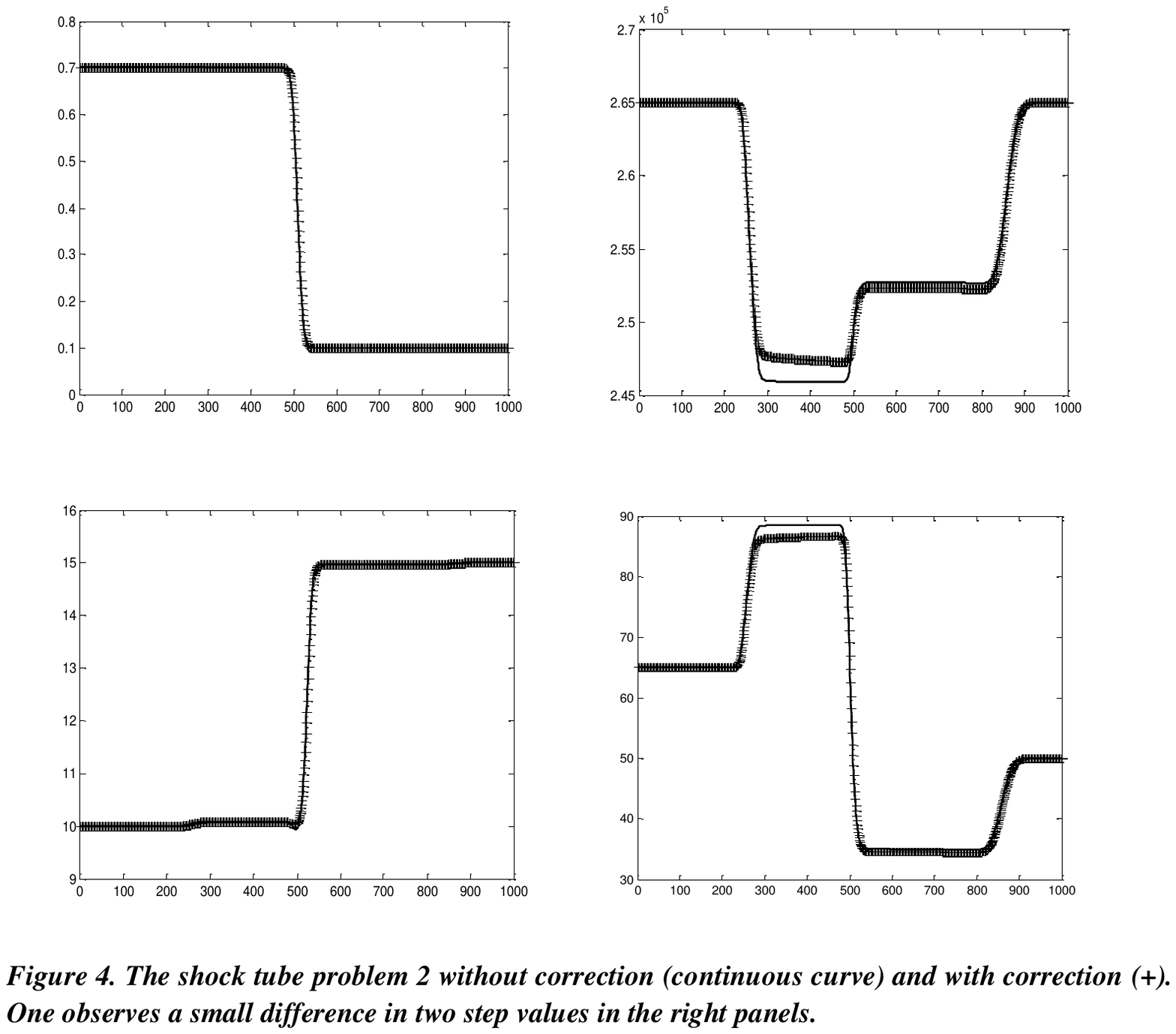}\end{figure}


\textit{Numerical observations.} First it has been observed that the scheme has always given the same result as the weak asymptotic method. It has the advantage  to be more  efficient and of a very easy use since one has only to fix the value of the CFL number $r$ and then the value of the averaging parameter $\mu$ in (52, 53).\\

The scheme in this section has been used with the interface pressure modelling (11) in \cite{EvjeFlatten} p. 180 which ensures the hyperbolicity of the system. In the case of shock tube problem 2 one can observe a slight difference relatively to the absence of correction (figure 4: 1000 space steps, $r=0.002, \mu=0.1$):
\vskip9cm
\textit{Figure 5. Quality of the transport-correction scheme: +++ results with pressure correction and (continuous line) without pressure correction. The curves are obtained with 100 space steps only.}\\
\\
 the second step values from the left in pressure and gas velocity are 248000 and 86.5 respectively instead of 246000 and 88.5. Since the tests in \cite{EvjeFlatten} have been done in presence of this pressure correction and are  close to the values we obtain with this correction, this explains the small disagreement observed when comparing the results in \cite{EvjeFlatten} figures 4 and 5 p. 197 and 198  with those in figure 2 for these two step values. Besides this difference the results in figure 1 and 3 are unchanged in absence or presence of the pressure correction, in particular the presence of the middle "singular  wave" is independent of the presence of pressure correction.\\

The numerical quality of the scheme is tested in figure 5, both in absence and presence of pressure correction: a dicretization in 100 space cells  suffices to obtain the step values and the jump formulas (as in the above arrays corresponding to figures 1, 2 and 3) with precision.\\

 \begin{figure}[h]\centering\includepdf[width=\textwidth]{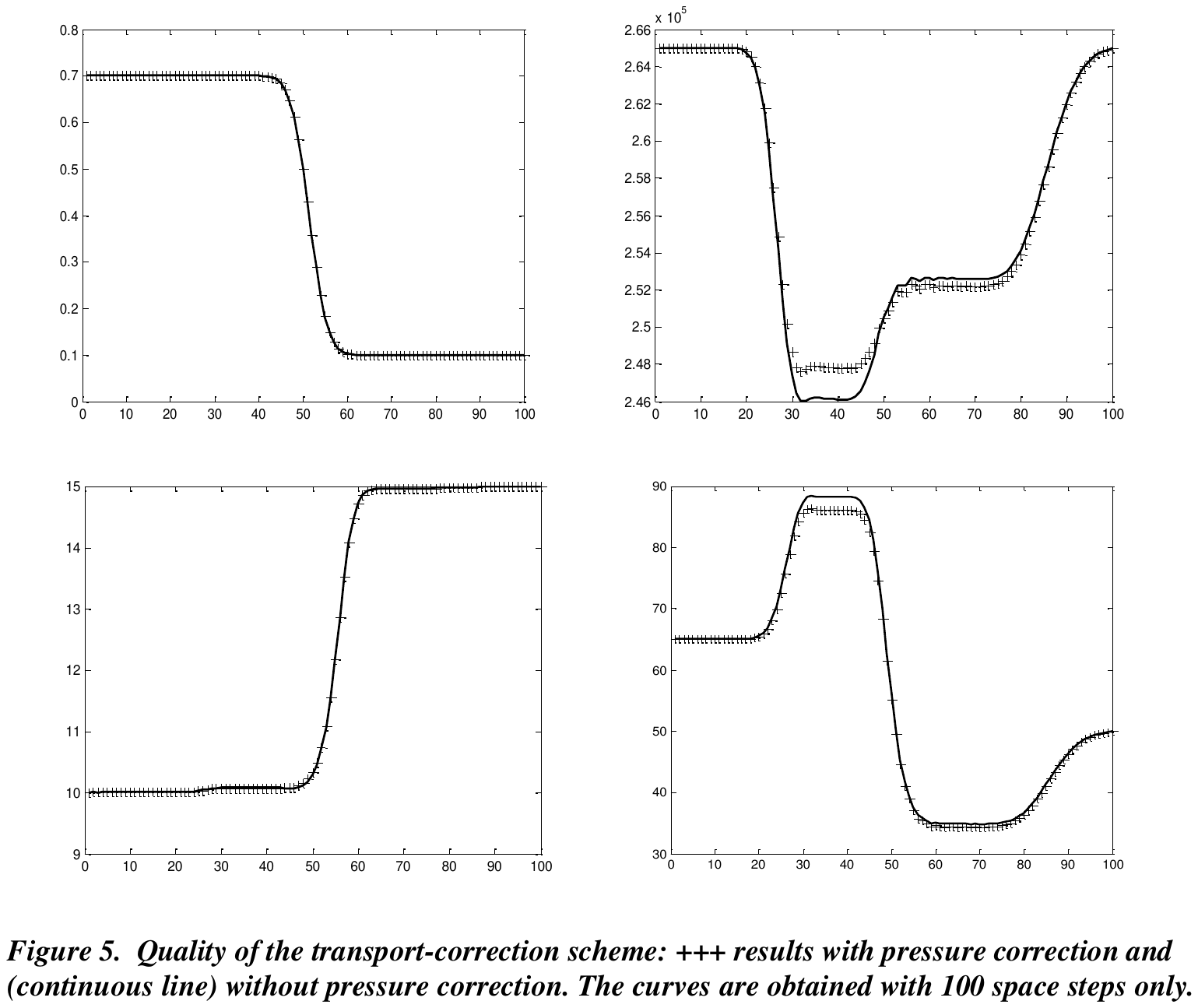}\end{figure}

\textbf{9. Conclusion.}\\
 The approximate solutions we have constructed with full proof and rather arbitrary initial data provide a mathematical tool that permits theorical and numerical investigations of the initial value problem for the equal pressure model of multifluid flows in the isothermal case. Since numerical calculations of these approximate solutions can be done easily and accurately with standard ODEs methods these approximate solutions can play the role of explicit solutions for mathematical  and  numerical investigations. They show that  numerical schemes from scientific computing give an approximate solution   besides the mathematical peculiarities of the model. They can show that supplementary terms such as pressure corrections  do not modify (shock tube problem 1) or modify only slightly (shock tube problem 2) the solution. They permit to investigate the nature of the "solutions" put in evidence by these approximate solutions and by scientific computing although this system is in nondivergence form.   \\
\\
 Acknowledgements. The author is very grateful to members of the Instituto de Matematica of The Universidade de S\~ao Paulo, of the Instituto de Matematica, Estatistica e Computa\c c\~ao Cientifica of the Universidade Estadual de Campinas and of the Instituto de Matematica of the Universidade Federal do Rio de Janeiro for their attention, encouragements and suggestions while doing this work.\\

\textbf{Appendix. Formal calculations on the system.}\\
 We obtain jump formulas from formal calculations. We observe in section 6 that these jump formulas are satisfied by the left shock waves and by the right shock waves in figures 1, 2 and 3. Developping (10) with $g=0$ and simplifying from (8), then dividing by $r_1$ one obtains\\
\begin{equation}\frac{\partial}{\partial t}(u_1)+u_1\frac{\partial}{\partial x}(u_1)+\frac{\alpha\frac{\partial}{\partial x} p}{r_1}=0.\end{equation}
Using the state law (6) $p=K_1\rho_1-b_1$ and $\rho_1=\frac{r_1}{\alpha}$ one obtains 
$$\frac{\partial}{\partial t} (u_1)=\frac{\partial}{\partial x}(-K_1log(\rho_1)-\frac{(u_1)^2}{2}),$$
which gives the jump condition
\begin{equation} c=K_1\frac{log(\rho_{1,r})-log(\rho_{1,l})}{[u_1]}+\frac{u_{1,r}+u_{1,l}}{2}
\end{equation}
where $c$ denotes the velocity of the shock wave. The same calculation holds from (10) and (9) and gives (62)  with index 2 and the same fomula with index 2.\\

\end{document}